\documentclass{tran-l}
\usepackage{amssymb}

\usepackage{amsmath}

\input{tcilatex}
\theoremstyle{definition}
\theoremstyle{remark}
\numberwithin{equation}{section}

\input tcilatex

\begin{document}
\title[A functional expression for the curvature...]{A functional expression for the curvature of hyper-dimensional \textit{%
Riemannian} spaces}
\author{\textbf{branko saric}}
\address{\textit{The Institute ''Kirilo Savic``, 11000 Belgrade, V. Stepe 51., Serbia}}
\email{bsaric@ptt.yu}
\date{April 26, 2000}
\maketitle

\begin{abstract}
Analogously to the concept of a curvature of curve and surface, in the
differential geometry, in the main part of this paper the concept of the
curvature of the hyper-dimensional vector spaces of \textit{Riemannian}
metric is generally defined. The defined concept of the curvature of \textit{%
Riemannian} spaces of higher dimensions $M$\textit{:}$\,M\geq 2$, in the
further text of the paper, is functional related to the fundamental
parameters of an internal geometry of space, more exactly, to components of 
\textit{Riemann-Christoffel's} tensor of curvature. At the end, analogously
to the concept of lines of curvature in the differential geometry, the
concept of sub-spaces of curvature of \textit{Riemannian} hyper-dimensional
vector spaces is also generally defined.
\end{abstract}

\section{Introduction}

The well-known \textit{Riemannian} mathematical model of defining the
curvature of hyper-dimensional curvilinear metric spaces, via \textit{%
Gaussian} concept of the two-dimensional surface curvature, \cite{An}\textit{%
;} \cite{L-C} and \cite{Pa}, is an imperfect in spite of that it functional
related to internal geometry of a space\textit{:}

1. Firstly, for a reason that the concept of the curvature of
hyper-dimensional metric spaces reduced to the concept of \textit{Gaussian}
curvature of two-dimensional geodesic surface, it cannot be said that 
\textit{Riemannian} concept of the curvature of hyper-dimensional spaces is
general one because it is an inapplicable to one-dimensional space.

2. As secondly, since there are more than one geodesic surface in vector
spaces of higher dimensions, it is obvious that \textit{Riemannian} concept
of the curvature of hyper-dimensional spaces, in the general case, is not
uniquely defined one.

3. Finally, at any an individual, concrete case, from the practical point of
view, it is not simple to come to the quantitatively usable functional
expression for \textit{Riemannian} curvature of an analyzed
hyper-dimensional metric space.

Hence, the one other mathematical model of defining the curvature of the
hyper-dimensional metric spaces, which essentially differs from \textit{%
Riemannian} model, is presented in this paper. Concretely, the mathematical
model is being discussed, which can be said to be generalization of
well-known model of defining the curvature of curve and surface in the
differential geometry.

\subsection{Basic characteristics of space continuum}

The concept of geometrical point is one of fundamental concepts. Closely
related to the concept of geometrical points is the system of values $%
(a_{1},a_{2},...\,a_{N})$ of some an arbitrary $N$ variables $%
(x_{1},x_{2},...\,x_{N})$ such that a set of all geometrical points, and for
all real values of the variables, is a real $N$ - dimensional space of a
space continuum, \cite{An}. The geometrical point $\mathbf{O}$, defined by
system of zero values $(0,0,...\,0)$, is an origin of system of reference
(of co-ordinate system) of the space. The vector $\vec{r}\left( x^{i}\right) 
$, defined with respect to the origin $\mathbf{O}$, is a position vector.
Note that the concept of a vector, in the vector hyper-dimensional spaces $%
\left( N>3\right) $, should be conditionally comprehend in the sense of its
geometrical presentation in a form of segments, hence it bears a name linear
tensor, \cite{Pa}.

Covariant vectors $\vec{e}_{i}$\textit{:} $\vec{e}_{i}=\partial _{x_{i}}\vec{%
r}\left( x^{j}\right) $, where $\partial _{x_{i}}=\frac{\partial }{\partial
x^{i}}$, form the covariant vector basis $\left\{ \vec{e}_{i}\right\}
_{i=1}^{N}$ of the tangent space of a space continuum. The vectors $\vec{e}%
^{i}$, such that at each point of a space\textit{:} $\vec{e}_{i}\cdot \vec{e}%
^{k}=\delta _{i}^{k}$, where the second order system of the unit values $%
\delta _{i}^{k}$ - is an unit $N\times N$ matrix (\textit{Kronecker's}
delta-symbol), \cite{An}\textit{;} \cite{L-C} and \cite{Pa}, form a basis $%
\left\{ \vec{e}^{i}\right\} _{i=1}^{N}$, which is called the dual basis of
the covariant vector basis $\left\{ \vec{e}_{i}\right\} _{i=1}^{N}$. This is
so-called natural isomorphism from $\left\{ \vec{e}_{i}\right\} _{i=1}^{N}$
onto $\left\{ \vec{e}^{i}\right\} _{i=1}^{N}$. The infinitesimal value $d%
\vec{r}$ of the position vector $\vec{r}$ of representative point is defined
by $d\vec{r}=dx^{i}\vec{e}_{i}=dx_{i}\vec{e}^{i}$, where the well-known 
\textit{Einstein's} convention is applied to summation with respect to the
repetitive indexes (uppers and lowers), \cite{An} and \cite{Pa}, herein as
well as in the further text of the paper.

\section{The main results}

\subsection{A curvature of hyper-dimensional \textit{Riemannian }spaces}

By the following transformation low\textit{:} $x^{i}=x^{i}\left( q^{\alpha
}\right) $\textit{;} $i=1,2,...,N$, $\alpha =1,2,...,M\leq N$, in the
general case, an arbitrary $M$ - dimensional metric space is defined, \cite
{An} 
\begin{equation*}
ds^{2}=\partial _{x^{i}}\vec{r}\partial _{q^{\alpha }}x^{i}\cdot \partial
_{x^{j}}\vec{r}\partial _{q^{\beta }}x^{j}dq^{\alpha }dq^{\beta
}=e_{ij}\partial _{q^{\alpha }}x^{i}\partial _{q^{\beta }}x^{j}dq^{\alpha
}dq^{\beta }=g_{\alpha \beta }dq^{\alpha }dq^{\beta },
\end{equation*}
where the positive definite symetric square matrices $e_{ij}$\textit{:} $%
e_{ij}=\vec{e}_{i}\cdot \vec{e}_{j}$ and $g_{\alpha \beta }$\textit{:} 
\begin{equation*}
g_{\alpha \beta }=\partial _{q^{\alpha }}\vec{r}\cdot \partial _{q^{\beta }}%
\vec{r}=\vec{g}_{\alpha }\cdot \vec{g}_{\beta }=\partial _{x^{i}}\vec{r}%
\partial _{q^{\alpha }}x^{i}\cdot \partial _{x^{j}}\vec{r}\partial
_{q^{\beta }}x^{j}=\vec{e}_{i}\cdot \vec{e}_{j}\partial _{q^{\alpha
}}x^{i}\partial _{q^{\beta }}x^{j},
\end{equation*}
of degree\textit{:} $N$ and $M$, respectively, are fundamental (metric)
tensors of an ambient $N$ - dimensional \textit{Euclidean} (more exactly 
\textit{Cartesian)} space $x^{i}$, as well as, in the general case of 
\textit{Riemannian} covering map ($M<N$), of an internal $M$ - dimensional 
\textit{Riemannian} curvilinear space $q^{\alpha }$ immersed in it. As it is
well-known if $M=N$ a smooth map $x^{i}\rightarrow q^{\alpha }$ is called an
isometric immersion.

The smallest possible dimensional difference $C$\textit{:} $C=N-M$, is said
to define a class of the \textit{Riemannian} vector space $q^{\alpha }$, 
\cite{An}. Hence, \textit{Riemannian} space $q^{\alpha }$\textit{:} $%
q^{\alpha }=q$, of class $C$\textit{:} $C=N-1$, is an arbitrary curve of the
ambient $N$ - dimensional \textit{Euclidean (Cartesian)} space $x^{i}$.

Since the vector $d_{q}\vec{g}$, where $d_{q}=\frac{d}{dq}$, as derivative
of the fundamental vector $\vec{g}$\textit{:} $\vec{g}=d_{q}\vec{r}$, lying
in the tangent space of \textit{Riemannian} space $q$ of class $C$\textit{:} 
$C=N-1$ (on the tangent of curve), along the curve, from the point of view
of the interior of an ambient space $x^{i}$, is a vector of the ambient
space $x^{i}$, then at each points of the curve there exist $C$\textit{:} $%
C=N-1$, the linearly independent and mutually orthogonal unit vectors\textit{%
:} $\vec{n}^{\Lambda }$\textit{;} $\Lambda =1,2,...,N-1$, being orthogonal
on the unit vector of tangent $\vec{t}$\textit{:} $\vec{t}=d_{s}\vec{r}$, of
the curve, such that 
\begin{equation}
d_{qq}^{2}\vec{r}=d_{q}\left( d_{q}x^{i}\vec{e}_{i}\right) =\partial _{x^{j}}%
\vec{e}_{i}d_{q}x^{j}d_{q}x^{i}+d_{qq}^{2}x^{i}\vec{e}_{i}=  \label{1}
\end{equation}
\begin{equation*}
=d_{q}\vec{g}=\left( d_{q}\vec{g}\cdot \vec{n}^{\Lambda }\right) \vec{n}%
_{\Lambda }+\left( d_{q}\vec{g}\cdot \vec{g}\right) \vec{g}.
\end{equation*}

The covariant vectors $\vec{n}_{\Sigma }$, as well as the vectors $\vec{n}%
^{\Lambda }$ satisfying the condition\textit{:} $\vec{n}_{\Sigma }\cdot \vec{%
n}^{\Lambda }=\delta _{\Sigma }^{\Lambda }$, form covariant $\left\{ \vec{n}%
_{\Sigma }\right\} _{\Sigma =1}^{N-1}$, as well as dual vector basis $%
\left\{ \vec{n}^{\Lambda }\right\} _{\Lambda =1}^{N-1}$ of a normal vector
space of \textit{Riemannian} space $q$ of class $C$\textit{:} $C=N-1$,
respectively.

From the point of view of the internal geometry of \textit{Riemannian} space 
$q$, the vector $\vec{\kappa}$\textit{:} 
\begin{equation*}
\vec{\kappa}=\left( d_{q}\vec{g}\cdot \vec{n}^{\Lambda }\right) \vec{n}%
_{\Lambda },
\end{equation*}
is a vector of curvature of curve, more exactly, from the point of view of
the internal geometry of an ambient \textit{Euclidean} space $x^{i}$, it is
a vector of geodesic curvature of curve, \cite{An}. Namely, by projection of
the vector $\vec{\kappa}$ onto the normal vector space of \textit{Riemannian}
space $q$ and onto the tangent vector space of ambient \textit{Euclidean}
space $x^{i}$, it is obtained that 
\begin{equation*}
\vec{\kappa}\cdot \vec{n}^{\Sigma }=d_{q}\vec{g}\cdot \vec{n}^{\Sigma
}=\left( \partial _{x^{j}}\vec{e}_{i}\cdot \vec{e}%
^{k}d_{q}x^{i}d_{q}x^{j}+d_{qq}^{2}x^{k}\right) \left( \vec{e}_{k}\cdot \vec{%
n}^{\Sigma }\right) =
\end{equation*}
\begin{equation*}
=\left( \Gamma _{ij}^{k}d_{q}x^{i}d_{q}x^{j}+d_{qq}^{2}x^{k}\right)
n_{k}^{\Sigma }
\end{equation*}
and 
\begin{equation*}
\vec{\kappa}\cdot \vec{e}^{l}=\left( d_{q}\vec{g}\cdot \vec{n}^{\Lambda
}\right) \vec{n}_{\Lambda }\cdot \vec{e}^{l}=\left( \Gamma
_{ij}^{k}d_{q}x^{i}d_{q}x^{j}+d_{qq}^{2}x^{k}\right) \vec{e}_{k}\cdot \vec{n}%
^{\Lambda }\left( \vec{n}_{\Lambda }\cdot \vec{e}^{l}\right) ,
\end{equation*}
where the mixed systems $\Gamma _{ij}^{k}$\textit{:} $\Gamma
_{ij}^{k}=\partial _{x^{j}}\vec{e}_{i}\cdot \vec{e}^{k}$, are coefficients
of connection (\textit{Christoffel's} symbols) of the second type, and an
intensity of vector 
\begin{equation*}
\frac{1}{\vec{g}\cdot \vec{g}}\left( \Gamma
_{ij}^{k}d_{q}x^{i}d_{q}x^{j}+d_{qq}^{2}x^{k}\right) \vec{e}_{k}=\left( d_{s}%
\vec{t}\cdot \vec{e}^{k}\right) \vec{e}_{k}=d_{s}\vec{t},
\end{equation*}
being in the direction of the vector $\vec{\kappa}$\textit{:} 
\begin{equation*}
\frac{1}{\vec{g}\cdot \vec{g}}\left( \Gamma
_{ij}^{k}d_{q}x^{i}d_{q}x^{j}+d_{qq}^{2}x^{k}\right) \left( \vec{e}_{k}\cdot 
\vec{n}^{\Lambda }\right) \vec{n}_{\Lambda }=\frac{1}{\vec{g}\cdot \vec{g}}%
\left( \vec{\kappa}\cdot \vec{e}^{l}\right) \vec{e}_{l}=\frac{1}{\vec{g}%
\cdot \vec{g}}\vec{\kappa},
\end{equation*}
defines both the curvature $\kappa $ of the curve, from the point of view of
the internal geometry of \textit{Riemannian} space $q$ of class $C$\textit{:}
$C=N-1$ and geodesic curvature $\kappa _{g}$ of the curve, from the point of
view of the internal geometry of ambient \textit{Euclidean} space $x^{i}$, 
\cite{An} and \cite{Pa}.

Analogously to the vector of curvature\textit{:} $\vec{\kappa}=\left( d_{q}%
\vec{g}\cdot \vec{n}^{\Lambda }\right) \vec{n}_{\Lambda }$, of a \textit{%
Riemannian} space $q$ of class $C$\textit{:} $C=N-1$, in the general case of 
\textit{Riemannian} space $q^{\alpha }$ of class $C$\textit{:} $C=N-M$,
immersed in $N$ - dimensional ambient \textit{Euclidean} space $x^{i}$, from
the viewpoint of its internal geometry, the matrix scheme of vectors 
\begin{equation}
\vec{K}_{\alpha \beta }=\left( \partial _{q^{\beta }}\vec{g}_{\alpha }\cdot 
\vec{n}^{\Lambda }\right) \vec{n}_{\Lambda }=\tau _{\alpha \beta }^{\Lambda }%
\vec{n}_{\Lambda },  \label{2}
\end{equation}
is a matrix scheme of the sectional curvature vectors, of a \textit{%
Riemannian} space $q^{\alpha }$ of class $C$\textit{:} $C=N-M$. The
determinant $\left| \tau _{\alpha \beta }^{\Lambda }\tau _{\gamma \delta
}^{\Sigma }g^{\delta \beta }n_{\Lambda \Sigma }\right| $ of the matrix
obtained by the matrix multiplication (rows $\times $ columns) of the matrix
schemes of curvature vectors\textit{:} $\tau _{\alpha \beta }^{\Lambda }\vec{%
n}_{\Lambda }$ and $\tau _{\gamma \delta }^{\Sigma }g^{\delta \beta }\vec{n}%
_{\Sigma }$, is the proportionate to the square of the curvature $\kappa $
of \textit{Riemannian} space $q^{\alpha }$ of class $C$\textit{:} $C=N-M$,
and with the coefficient of proportionality which is equal to the
determinant of matrix of fundamental tensor $g_{\alpha \beta }$%
\begin{equation*}
\left| \vec{K}_{\alpha \beta }\cdot \vec{K}_{\gamma }^{\beta }\right|
=\left| \tau _{\alpha \beta }^{\Lambda }\tau _{\gamma \delta }^{\Sigma
}g^{\delta \beta }n_{\Lambda \Sigma }\right| =\left| g_{\alpha \beta
}\right| \kappa ^{2},
\end{equation*}
more exactly, 
\begin{equation}
\kappa ^{2}=\frac{\left| \tau _{\alpha \beta }^{\Lambda }\tau _{\gamma
\delta }^{\Sigma }g^{\delta \beta }n_{\Lambda \Sigma }\right| }{\left|
g_{\alpha \beta }\right| }.  \label{3}
\end{equation}

\begin{description}
\item[Comment]  If $N$ - dimensional ambient space $x^{i}$ is either \textit{%
Riemannian} curvilinear space of class $C$ or \textit{Euclidean} curvilinear
space of class $C$\textit{:} $C\geq 1$, both with normal vector space $\vec{w%
}_{P}$\textit{;} $P=1,2,...,C$, then the matrix scheme of vectors 
\begin{equation*}
\vec{K}_{\alpha \beta }=\left( \partial _{q^{\beta }}\vec{g}_{\alpha }\cdot 
\vec{n}^{\Lambda }\right) \vec{n}_{\Lambda }+\left( \partial _{q^{\beta }}%
\vec{g}_{\alpha }\cdot \vec{w}^{P}\right) \vec{w}_{P},
\end{equation*}
is a matrix scheme of vectors of curvature of \textit{Riemannian} space $%
q^{\alpha }$ of class $\hat{C}$\textit{: }$\hat{C}=C+N-M$.

From the point of view of the internal geometry of the ambient space $x^{i}$%
, the partial matrix schemes of vectors of curvature\textit{:} $\vec{N}%
_{\alpha \beta }=\left( \partial _{q^{\beta }}\vec{g}_{\alpha }\cdot \vec{w}%
^{P}\right) \vec{w}_{P}=\varkappa _{\alpha \beta }^{P}\vec{w}_{P}$ and $\vec{%
G}_{\alpha \beta }=\left( \partial _{q^{\beta }}\vec{g}_{\alpha }\cdot \vec{n%
}^{\Lambda }\right) \vec{n}_{\Lambda }=\tau _{\alpha \beta }^{\Lambda }\vec{n%
}_{\Lambda }$, are matrix schemes of vectors of normal curvature, as well as
of geodesic curvature, of \textit{Riemannian} space $q^{\alpha }$ of class $%
\hat{C}$, respectively.

The square of intensity of curvature $\kappa $ of \textit{Riemannian} space
of class $\hat{C}$, in this case is equal to the sum of square of the normal
curvature $\kappa _{n}$\textit{:} 
\begin{equation*}
\kappa _{n}=\frac{\left| \varkappa _{\alpha \beta }^{P}\varkappa _{\gamma
\delta }^{L}g^{\beta \delta }w_{PL}\right| }{\left| g^{\alpha \beta }\right| 
}
\end{equation*}
and of the geodesic curvature $\kappa _{g}$\textit{:} 
\begin{equation*}
\kappa _{g}=\frac{\left| \tau _{\alpha \beta }^{\Lambda }\tau _{\gamma
\delta }^{\Sigma }g^{\beta \delta }n_{\Lambda \Sigma }\right| }{\left|
g^{\alpha \beta }\right| },
\end{equation*}
of \textit{Riemannian} space $q^{\alpha }$ of class $\hat{C}$\textit{:} 
\begin{equation*}
\kappa ^{2}=\kappa _{n}^{2}+\kappa _{g}^{2}.
\end{equation*}

The space $q^{\alpha }$ of class $\hat{C}$, as well as of the geodesic
curvature zero $\kappa _{g}$\textit{:} $\kappa _{g}=0$, with respect to the
ambient space $x^{i}$ of class $C$, is a geodesic sub-space of that ambient
space.$\blacktriangledown $
\end{description}

\subsubsection{Functional expression for the curvature of hyper-dimensional 
\textit{Riemannian} spaces, from the point of view of the internal geometry
of space.}

If the vector functional equality (\ref{1}) 
\begin{equation*}
\partial _{q^{\beta }}\vec{g}_{\alpha }=\left( \partial _{q^{\beta }}\vec{g}%
_{\alpha }\cdot \vec{g}^{\gamma }\right) \vec{g}_{\gamma }+\left( \partial
_{q^{\beta }}\vec{g}_{\alpha }\cdot \vec{n}^{\Lambda }\right) \vec{n}%
_{\Lambda },
\end{equation*}
in the general case of the ambient \textit{Riemannian} space $q^{\alpha }$
of class $C$\textit{:} $C=N-M$, is partially differentiated with respect to
dual co-ordinates $q^{\delta }$, it is obtained that 
\begin{equation}
\partial _{q^{\beta }q^{\delta }}^{2}\vec{g}_{\alpha }=\partial _{q^{\delta
}}\Gamma _{\alpha \beta }^{\gamma }\vec{g}_{\gamma }+\Gamma _{\alpha \beta
}^{\gamma }\partial _{q^{\delta }}\vec{g}_{\gamma }+\partial _{q^{\delta
}}\tau _{\alpha \beta }^{\Lambda }\vec{n}_{\Lambda }+\tau _{\alpha \beta
}^{\Lambda }\partial _{q^{\delta }}\vec{n}_{\Lambda },  \label{4}
\end{equation}
that is, 
\begin{equation}
\partial _{q^{\delta }q^{\beta }}^{2}\vec{g}_{\alpha }=\partial _{q^{\beta
}}\Gamma _{\alpha \delta }^{\gamma }\vec{g}_{\gamma }+\Gamma _{\alpha \delta
}^{\gamma }\partial _{q^{\beta }}\vec{g}_{\gamma }+\partial _{q^{\beta
}}\tau _{\alpha \delta }^{\Lambda }\vec{n}_{\Lambda }+\tau _{\alpha \delta
}^{\Lambda }\partial _{q^{\beta }}\vec{n}_{\Lambda }.  \label{5}
\end{equation}

Since differentials $d\vec{g}_{\alpha }=\partial _{q^{\beta }}\vec{g}%
_{\alpha }dq^{\beta }$ are an absolute ones, more exactly, the covariant
vector basis $\left\{ \vec{g}_{\alpha }\right\} $ of a \textit{Riemannian}
space $q^{\alpha }$ of class $C$\textit{:} $C=N-M$, is uniquely defined at
each point of a space, and in accordance with that the condition of
integrability\textit{:} $\partial _{q^{\beta }q^{\delta }}^{2}\vec{g}%
_{\alpha }-\partial _{q^{\delta }q^{\beta }}^{2}\vec{g}_{\alpha }=0$, is
satisfied, then by projection of difference of the respectable vector
functional expressions on both sides of the previous equations\textit{:} (%
\ref{4}) and (\ref{5}), onto covariant vector basis $\left\{ \vec{g}_{\rho
}\right\} $, the functional relation is obtained 
\begin{equation*}
\partial _{q^{\delta }}\Gamma _{\alpha \beta }^{\gamma }g_{\gamma \rho
}-\partial _{q^{\beta }}\Gamma _{\alpha \delta }^{\gamma }g_{\gamma \rho
}+\Gamma _{\alpha \beta }^{\lambda }\Gamma _{\lambda \delta }^{\gamma
}g_{\gamma \rho }-\Gamma _{\alpha \delta }^{\lambda }\Gamma _{\lambda \beta
}^{\gamma }g_{\gamma \rho }=\tau _{\alpha \beta }^{\Lambda }\tau _{\delta
\rho ,\Lambda }-\tau _{\alpha \delta }^{\Lambda }\tau _{\beta \rho ,\Lambda
},
\end{equation*}
having in mind both the equation (\ref{2}) and orthogonality of vectors%
\textit{:} $\vec{g}_{\alpha }$ and $\vec{n}_{\Lambda }$. On account of the
fact that the functional expression on the left hand side of preceding
relation represents \textit{Riemann-Christoffel's} tensor of \textit{%
Riemannian} space $q^{\alpha }$ of class $C$\textit{:} $C=N-M$, \cite{An}%
\textit{;} \cite{L-C} and \cite{Pa}, it finally follows that 
\begin{equation}
R_{\rho \alpha \delta \beta }=\tau _{\alpha \beta }^{\Lambda }\tau _{\delta
\rho ,\Lambda }-\tau _{\alpha \delta }^{\Lambda }\tau _{\beta \rho ,\Lambda
}.  \label{6}
\end{equation}

By application of \textit{Gauss-Chi\`{o}'s} procedure for condensation of
determinants, \cite{M-D}, and on the basis of tensorial relation (\ref{6}),
as well as of functional relation (\ref{3}), it is possible the curvature of 
\textit{Riemannian} space $q^{\alpha }$ of class $C$\textit{:} $C=N-M$, to
come into functional relation to internal geometry of a space, more exactly
to the components of \textit{Rimann-Christoffel's} tensor of the curvature $%
R_{\rho \alpha \delta \beta }$.

Namely, if in addition to the matrix scheme of vectors\textit{:} $\vec{K}%
_{\alpha \beta }=\tau _{\alpha \beta }^{\Lambda }\vec{n}_{\Lambda }$, the
matrix schemes of vectors\textit{:} 
\begin{equation*}
\vec{\Phi}_{\alpha \beta }=\left[ 
\begin{array}{lllll}
\phi ^{\Lambda }\vec{n}_{\Lambda } & 0 & ... & 0 & 0 \\ 
\tau _{12}^{\Lambda }\vec{n}_{\Lambda } & \tau _{11}^{\Lambda }\vec{n}%
_{\Lambda } & ... & 0 & 0 \\ 
... & ... & ... & ... & ... \\ 
\tau _{1\left( M-1\right) }^{\Lambda }\vec{n}_{\Lambda } & 0 & ... & \tau
_{11}^{\Lambda }\vec{n}_{\Lambda } & 0 \\ 
\tau _{1M}^{\Lambda }\vec{n}_{\Lambda } & 0 & ... & 0 & \tau _{11}^{\Lambda }%
\vec{n}_{\Lambda }
\end{array}
\right] ;
\end{equation*}
\begin{equation*}
\vec{Z}_{\alpha \beta }=\left[ 
\begin{array}{lllll}
\tau _{12}^{\Lambda }\vec{n}_{\Lambda } & 0 & ... & 0 & 0 \\ 
\tau _{11}^{\Lambda }\vec{n}_{\Lambda } & \zeta ^{\Lambda }\vec{n}_{\Lambda }
& ... & \tau _{1\left( M-1\right) }^{\Lambda }\vec{n}_{\Lambda } & \tau
_{1M}^{\Lambda }\vec{n}_{\Lambda } \\ 
... & ... & ... & ... & ... \\ 
0 & 0 & ... & \tau _{12}^{\Lambda }\vec{n}_{\Lambda } & 0 \\ 
0 & 0 & ... & 0 & \tau _{12}^{\Lambda }\vec{n}_{\Lambda }
\end{array}
\right] ;
\end{equation*}
\begin{equation*}
\vec{\Xi}_{\alpha \beta }=\left[ 
\begin{array}{lllll}
\xi ^{\Lambda }\vec{n}_{\Lambda } & 0 & ... & 0 & 0 \\ 
\tau _{22}^{\Lambda }\vec{n}_{\Lambda } & \tau _{21}^{\Lambda }\vec{n}%
_{\Lambda } & ... & 0 & 0 \\ 
... & ... & ... & ... & ... \\ 
\tau _{2\left( M-1\right) }^{\Lambda }\vec{n}_{\Lambda } & 0 & ... & \tau
_{21}^{\Lambda }\vec{n}_{\Lambda } & 0 \\ 
\tau _{2M}^{\Lambda }\vec{n}_{\Lambda } & 0 & ... & 0 & \tau _{21}^{\Lambda }%
\vec{n}_{\Lambda }
\end{array}
\right] ;
\end{equation*}
\begin{equation*}
\vec{\Psi}_{\alpha \beta }=\left[ 
\begin{array}{lllll}
\tau _{22}^{\Lambda }\vec{n}_{\Lambda } & 0 & ... & 0 & 0 \\ 
\tau _{21}^{\Lambda }\vec{n}_{\Lambda } & \psi ^{\Lambda }\vec{n}_{\Lambda }
& ... & \tau _{2\left( M-1\right) }^{\Lambda }\vec{n}_{\Lambda } & \tau
_{2M}^{\Lambda }\vec{n}_{\Lambda } \\ 
... & ... & ... & ... & ... \\ 
0 & 0 & ... & \tau _{22}^{\Lambda }\vec{n}_{\Lambda } & 0 \\ 
0 & 0 & ... & 0 & \tau _{22}^{\Lambda }\vec{n}_{\Lambda }
\end{array}
\right] ,
\end{equation*}
are also introduced into analysis, where\textit{:} $\phi ^{\Lambda }\vec{n}%
_{\Lambda }$\textit{;}$\zeta ^{\Lambda }\vec{n}_{\Lambda }$\textit{;}$\xi
^{\Lambda }\vec{n}_{\Lambda }$ and $\psi ^{\Lambda }\vec{n}_{\Lambda }$, are
arbitrary vector functions of the normal vector space $\vec{n}_{\Lambda }$
of \textit{Riemannian} space $q^{\alpha }$ of class $C$\textit{:} $C=N-M$,
then having in view the fact that determinant of the product of matrices, as
it is well-known, is equal to the product of determinants of any matrix
separately, \cite{An} and \cite{M-D}, it follows that 
\begin{equation}
\left( \vec{\Phi}_{\gamma \eta }\cdot \vec{K}_{\delta }^{\eta }\right)
\left( \vec{K}_{\alpha \lambda }\cdot \vec{\Psi}_{\beta }^{\lambda }\right) =%
\frac{\left( \tau _{11}^{\Lambda }\phi _{\Lambda }\right) \left( \tau
_{22}^{\Lambda }\psi _{\Lambda }\right) }{\left| g_{\alpha \beta }\right|
^{2}}\times  \label{7}
\end{equation}
\begin{equation*}
\times \left| 
\begin{array}{llll}
R_{1212} & R_{1213} & ... & R_{121M} \\ 
R_{1312} & R_{1313} & ... & R_{131M} \\ 
... & ... & ... & ... \\ 
R_{1M12} & R_{1M13} & ... & R_{1M1M}
\end{array}
\right| \left| 
\begin{array}{llll}
R_{1212} & R_{1232} & ... & R_{12M2} \\ 
R_{3212} & R_{3232} & ... & R_{32M2} \\ 
... & ... & ... & ... \\ 
R_{M212} & R_{M232} & ... & R_{M2M2}
\end{array}
\right| =
\end{equation*}
\begin{equation*}
=\frac{\left( \tau _{11}^{\Lambda }\phi _{\Lambda }\right) \left( \tau
_{22}^{\Lambda }\psi _{\Lambda }\right) }{\left| g_{\alpha \beta }\right|
^{2}}\Delta _{\tau _{11}^{\Lambda }}^{R}\Delta _{\tau _{22}^{\Lambda }}^{R},
\end{equation*}
more exactly\footnote{%
\begin{equation*}
\left| \vec{\Phi}_{\alpha \delta }\cdot \vec{\Psi}_{\beta }^{\delta }\right|
=\frac{1}{\left| g_{\alpha \beta }\right| }\left| 
\begin{array}{llll}
\phi ^{\Lambda }\tau _{22,\Lambda } & 0 & ... & 0 \\ 
-\left( \tau _{22}^{\Lambda }\tau _{12,\Lambda }+\tau _{11}^{\Lambda }\tau
_{21,\Lambda }\right) & \psi ^{\Lambda }\tau _{11,\Lambda } & ... & -\tau
_{11}^{\Lambda }\tau _{2M,\Lambda } \\ 
... & ... & ... & ... \\ 
-\tau _{22}^{\Lambda }\tau _{1M,\Lambda } & 0 & ... & \tau _{11}^{\Lambda
}\tau _{22,\Lambda }
\end{array}
\right| =
\end{equation*}
\begin{equation*}
=\frac{1}{\left| g_{\alpha \beta }\right| }\left( \phi ^{\Lambda }\tau
_{22,\Lambda }\right) \left( \psi ^{\Lambda }\tau _{11,\Lambda }\right)
\left( \tau _{11}^{\Lambda }\tau _{22,\Lambda }\right) ^{M-2}.
\end{equation*}
} 
\begin{equation*}
\left| \vec{K}_{\alpha \delta }\cdot \vec{K}_{\beta }^{\delta }\right|
\left( \tau _{11}^{\Lambda }\tau _{22,\Lambda }\right) ^{M-2}=\frac{1}{%
\left| g_{\alpha \beta }\right| }\left| 
\begin{array}{llll}
R_{1212} & R_{1213} & ... & R_{121M} \\ 
R_{1312} & R_{1313} & ... & R_{131M} \\ 
... & ... & ... & ... \\ 
R_{1M12} & R_{1M13} & ... & R_{1M1M}
\end{array}
\right| \times
\end{equation*}
\begin{equation*}
\times \left| 
\begin{array}{llll}
R_{1212} & R_{1232} & ... & R_{12M2} \\ 
R_{3212} & R_{3232} & ... & R_{32M2} \\ 
... & ... & ... & ... \\ 
R_{M212} & R_{M232} & ... & R_{M2M2}
\end{array}
\right| =\frac{1}{\left| g_{\alpha \beta }\right| }\Delta _{\tau
_{11}^{\Lambda }}^{R}\Delta _{\tau _{22}^{\Lambda }}^{R},
\end{equation*}
considering the fact that the functional expression in the relation (\ref{7}%
), which is an independent of the choice of arbitrary vector functions%
\textit{:} $\phi ^{\Lambda }\vec{n}_{\Lambda }$ and $\psi ^{\Lambda }\vec{n}%
_{\Lambda }$ (the specially interesting case is one in which $\phi ^{\Lambda
}=\psi ^{\Lambda }=\tau _{12}^{\Lambda }$), defines the determinant $\left| 
\vec{K}_{\alpha \delta }\cdot \vec{K}_{\beta }^{\delta }\right| $.

Similarly 
\begin{equation*}
\left| \vec{K}_{\alpha \delta }\cdot \vec{K}_{\beta }^{\delta }\right|
\left( \tau _{12}^{\Lambda }\tau _{21,\Lambda }\right) ^{M-2}=\frac{1}{%
\left| g_{\alpha \beta }\right| }\left| 
\begin{array}{llll}
R_{1212} & R_{1232} & ... & R_{12M2} \\ 
R_{1213} & R_{1233} & ... & R_{12M3} \\ 
... & ... & ... & ... \\ 
R_{121M} & R_{123M} & ... & R_{12MM}
\end{array}
\right| \times
\end{equation*}
\begin{equation*}
\times \left| 
\begin{array}{llll}
R_{2121} & R_{2123} & ... & R_{212M} \\ 
R_{2131} & R_{2133} & ... & R_{213M} \\ 
... & ... & ... & ... \\ 
R_{21M1} & R_{21M3} & ... & R_{21MM}
\end{array}
\right| =\frac{1}{\left| g_{\alpha \beta }\right| }\Delta _{\tau
_{12}^{\Lambda }}^{R}\Delta _{\tau _{21}^{\Lambda }}^{R}.
\end{equation*}

On the basis of the functional formulation of the curvature of \textit{%
Riemannian} space $q^{\alpha }$ of class $C$\textit{:} $C=N-M$, the relation
(\ref{3})\textit{:} $\kappa ^{2}=\frac{\left| \vec{K}_{\alpha \delta }\cdot 
\vec{K}_{\beta }^{\delta }\right| }{\left| g_{\alpha \beta }\right| }$, as
well as of the two previously derived relations, and for $\left( M>2\right) $%
, it finally follows that 
\begin{equation}
\kappa ^{2}=\frac{1}{\left| g_{\alpha \beta }\right| ^{2}\left(
R_{1212}\right) ^{M-2}}\left[ \left( \Delta _{\tau _{11}^{\Lambda
}}^{R}\Delta _{\tau _{22}^{\Lambda }}^{R}\right) ^{\frac{1}{M-2}}-\left(
\Delta _{\tau _{12}^{\Lambda }}^{R}\Delta _{\tau _{21}^{\Lambda
}}^{R}\right) ^{\frac{1}{M-2}}\right] ^{M-2}.  \label{8}
\end{equation}

Clearly, in the case of \textit{Riemannian} space $q^{\alpha }$ of class $C$%
\textit{:} $C=N-M$ ($M=2$), functional expression (\ref{3}) for the
curvature of surface is reduced to the well-known \textit{Gaussian}
curvature of surface in the theory of surfaces, \cite{An} and \cite{Vi} 
\begin{equation*}
\kappa =\frac{R_{1212}}{\left| g_{\alpha \beta }\right| }.
\end{equation*}

\begin{description}
\item[Comment]  In the case when the component of \textit{%
Riemann-Christoffel's} tensor is equal to zero, it is possible, in the
functional expression (\ref{8}) for the curvature of \textit{Riemannian}
spaces, to take any another combination of the components of the matrix
scheme of the curvature vectors of space\textit{:} $\vec{K}_{\alpha \beta
}=\tau _{\alpha \beta }^{\Lambda }\vec{n}_{\Lambda }$, as the support
elements of \textit{Gauss-Chi\`{o}'s} procedure for condensation of
determinants, and for which some of the components of \textit{%
Riemann-Christoffel's} tensor of curvature are not equal to zero. Clearly,
if all components of \textit{Riemann-Christoffel's} tensor of the curvature
of space are identically equal to zero, then the space is \textit{Riemannian}
space of curvature zero (\textit{Euclidean} space).$\blacktriangledown $
\end{description}

\subsection{ Sub-spaces of curvature of hyper-dimensional spaces}

Let $N$ - dimensional ambient space of space continuum $x^{i}$, just as in
the first \textit{Comment} of preceding \textit{Section 2.1} of this paper,
be either \textit{Riemannian} space of class $C$ or \textit{Euclidean}
curvilinear space of class $C$\textit{:} $C\geq 1$, both with the normal
vector space $\vec{w}_{P}$\textit{;} $P=1,2,...,C$.

If and only if the vector components\textit{:} 
\begin{equation*}
\left( \partial _{q^{\delta }}\vec{N}_{\alpha \beta }g^{\alpha \delta }\cdot 
\vec{e}_{i}\right) \vec{e}^{i}=\varkappa _{\alpha \beta }^{P}g^{\alpha
\delta }\left( \partial _{q^{\delta }}\vec{w}_{P}\cdot \vec{e}_{i}\right) 
\vec{e}^{i},
\end{equation*}
of vectors obtained by a partial differentiation of the matrix scheme of
vector of normal curvature\textit{:} $\vec{N}_{\alpha \beta }=\left(
\partial _{q^{\delta }}\vec{g}_{\alpha }\cdot \vec{w}^{P}\right) \vec{w}%
_{P}=\varkappa _{\alpha \beta }^{P}\vec{w}_{P}$, of \textit{Riemannian}
space $q^{\alpha }$ of class $\hat{C}$\textit{:} $\hat{C}=C+N-M$, immersed
in ambient space $x^{i}$ of class $C$, are vectors of the tangent vector
space $\vec{g}_{\alpha }$ of \textit{Riemannian} space $q^{\alpha }$ of
class $\hat{C}$, in other words if and only if the mutually equivalent
conditions\textit{:} 
\begin{equation}
-\left( \partial _{q^{\delta }}\vec{N}_{\alpha \beta }g^{\alpha \delta
}\cdot \vec{e}_{i}\right) \vec{e}^{i}=\hat{\kappa}^{2}e_{ij}\partial
_{q^{\beta }}x^{j}\vec{e}^{i}=\hat{\kappa}^{2}\left( \vec{g}_{\beta }\cdot 
\vec{e}_{i}\right) \vec{e}^{i}=\hat{\kappa}^{2}\vec{g}_{\beta }  \label{9}
\end{equation}
and 
\begin{equation}
-\varkappa _{\alpha \beta }^{P}g^{\alpha \delta }\left( \partial _{q^{\delta
}}\vec{w}_{P}\cdot \vec{e}_{i}\right) =\varkappa _{\alpha \beta
}^{P}g^{\alpha \delta }t_{ij,P}\partial _{q^{\delta }}x^{j}=\hat{\kappa}%
^{2}e_{ij}\partial _{q^{\beta }}x^{j},  \label{10}
\end{equation}
are satisfied\footnote{{\footnotesize The second of the two relations of
equality is reduced to: }$t_{k,P}^{j}t_{ij}^{P}\partial _{q^{\beta }}x^{k}=%
\hat{\kappa}^{2}e_{ij}\partial _{q^{\beta }}x^{j}${\footnotesize , clearly
on the condition that all vectors of vector components }$\left( \partial
_{q^{\delta }}\vec{w}^{P}\cdot \vec{e}^{i}\right) \vec{e}_{i}${\footnotesize %
\ are also vectors of the tangent vector space }$\vec{g}_{\alpha }$%
{\footnotesize . Namely, since: }$t_{k,P}^{j}\partial _{q^{\beta
}}x^{k}=\partial _{q^{\beta }}\vec{w}_{P}\cdot \vec{e}^{i}${\footnotesize \
and }$\varkappa _{\alpha \beta ,P}g^{\alpha \delta }\partial _{q^{\delta
}}x^{i}=\left( \partial _{q^{\beta }}\vec{w}_{P}\cdot \vec{g}_{\alpha
}\right) \vec{g}^{\alpha }\cdot \vec{g}^{\delta }\left( \vec{e}^{i}\cdot 
\vec{g}_{\delta }\right) ${\footnotesize , then }$\varkappa _{\alpha \beta
,P}g^{\alpha \delta }\partial _{q^{\delta }}x^{i}=t_{k,P}^{j}\partial
_{q^{\beta }}x^{k}${\footnotesize .}}, then \textit{Riemannian} sub-space $%
q^{\alpha }$ of class $\hat{C}$, is a sub-space of curvature of the ambient
space $x^{i}$ of class $C$. In the case in which a \textit{Riemannian} space 
$q^{\alpha }$ of class $\hat{C}$\textit{:} $\hat{C}=C+N-M$, is an arbitrary
curve $\left( M=1\right) $, the preceding conditions are reduced to the
conditions\textit{:} $-\left( \partial _{q}\vec{N}\frac{1}{\vec{g}\cdot \vec{%
g}}\cdot \vec{e}_{i}\right) \vec{e}^{i}=\hat{\kappa}^{2}e_{ij}\partial
_{q}x^{j}\vec{e}^{i}$ and $\left( d_{q}\vec{g}\frac{1}{\vec{g}\cdot \vec{g}}%
\cdot \vec{w}^{P}\right) t_{ij,P}\partial _{q}x^{j}\vec{e}^{i}=\hat{\kappa}%
^{2}e_{ij}\partial _{q}x^{j}\vec{e}^{i}$, as well as $\left( d_{s}\vec{t}%
\cdot \vec{w}^{P}\right) t_{ij,P}\partial _{s}x^{j}\vec{e}%
^{i}=k_{ij}\partial _{s}x^{j}\vec{e}^{i}=\hat{\kappa}^{2}e_{ij}\partial
_{s}x^{j}\vec{e}^{i}$.

The last of them is the well-known condition for the curve to be a line of
the curvature of the ambient space $x^{i}$ of class $C$, \cite{Vi}. It is
obvious from this that conditions\textit{:} (\ref{9}) and (\ref{10}), are a
generalization of the preceding conditions.

By projection of the condition (\ref{9}), onto covariant vector basis $%
\left\{ \vec{g}_{\gamma }\right\} $, it is obtained that 
\begin{equation*}
\varkappa _{\alpha \beta }^{P}g^{\alpha \delta }t_{ij,P}\partial _{q^{\delta
}}x^{j}\partial _{q^{\gamma }}x^{i}=\varkappa _{\alpha \beta }^{P}\varkappa
_{\gamma ,P}^{\alpha }=\hat{\kappa}^{2}e_{ij}\partial _{q^{\beta
}}x^{j}\partial _{q^{\gamma }}x^{i}=\hat{\kappa}^{2}g_{\beta \gamma },
\end{equation*}
more exactly, 
\begin{equation*}
\left( \kappa _{n}\right) ^{2}=\frac{\left| \varkappa _{\alpha \beta
}^{P}\varkappa _{\gamma ,P}^{\alpha }\right| }{\left| g_{\beta \gamma
}\right| }=\hat{\kappa}^{2M},
\end{equation*}
with regard to the relation (\ref{3}).

\begin{description}
\item[Comment]  \textit{Riemannian} space of curvature\textit{:} $q^{\alpha
} $, of class $\hat{C}$, as a sub-space of the ambient space $x^{i}$ of
class $C$, is said to be the principal, if and only if 
\begin{equation}
t_{j,P}^{k}t_{ik}^{P}\partial _{q^{\beta }}x^{j}=\hat{\kappa}%
^{2}e_{ij}\partial _{q^{\beta }}x^{j},  \label{11}
\end{equation}
in other words, if all vectors of the vector components\textit{:} $\left(
\partial _{q^{\delta }}\vec{w}_{P}\cdot \vec{e}^{i}\right) \vec{e}_{i}$, are
vectors of the tangent vector space (see the \textit{Footnote 2} of the
paper).

In order to exist nontrivial solutions of the previous homogeneous linear
system (\ref{11}) of ordinary differential equations with respect to the
unknowns $\partial _{q^{\beta }}x^{j}$, the following condition 
\begin{equation}
\left| t_{j,P}^{k}t_{ik}^{P}-\hat{\kappa}^{2}e_{ij}\right| =0,  \label{12}
\end{equation}
must be satisfied, \cite{An} and \cite{M-D}.

On the basis of the developed form of the preceding condition (\ref{12}),
expressed by polynomial of $N$-th degree with respect to the unknown $\hat{%
\kappa}^{2}$, as well as of \textit{Vi\`{e}te's} formulas \cite{M-D}, and
taking the relation (\ref{3}) into consideration, it finally follows that 
\begin{equation*}
k^{2}=\frac{\left| t_{j,P}^{k}t_{ik}^{P}\right| }{\left| e_{ij}\right| }=%
\underset{i=1}{\overset{N}{\prod }}\hat{\kappa}_{i}^{2}.
\end{equation*}

If the varied of all eigevalues $\hat{\kappa}_{i}$ of the matrix\textit{:} $%
t_{j,P}^{k}t_{ik}^{P}$, is an unit, then the principal \textit{Riemannian}
sub-spaces of curvature\textit{:} $q^{\alpha }$, of class $\hat{C}$, are
one-dimensional ones $\left( M=1\right) $, in other words there exist $N$
principal directions of curvature with the normal curvatures\textit{:} $%
\left( \kappa _{n}\right) _{i}=\hat{\kappa}_{i}$, of the ambient space $%
x_{i} $ of class $C$. On the condition that there exists at least one of all
eigevalues $\tilde{\kappa}$ with the varied $M$, then there exists at least
one principal $M$ - dimensional \textit{Riemannian} sub-space of curvature%
\textit{:} $q^{\alpha }$, of class $\hat{C}$, with the normal curvature%
\textit{:} $\kappa _{n}=\tilde{\kappa}^{M}$, as well as $N-M$ mutually
orthogonal principal directions of curvature with the normal curvatures%
\textit{:} $\left( \kappa _{n}\right) _{i}=\hat{\kappa}_{i}$ , such that are
orthogonal onto the principal $M$ - dimensional \textit{Riemannian}
sub-spaces of curvature\textit{:} $q^{\alpha }$, of class $\hat{C}$. In this
emphasized case 
\begin{equation*}
k=\tilde{\kappa}^{M}\underset{i=1}{\overset{N-M}{\prod }}\hat{\kappa}_{i}.
\end{equation*}

Furthermore, if the vector components\textit{:} 
\begin{equation*}
\left( \partial _{q^{\delta }}\vec{N}_{\alpha \beta }\cdot \vec{e}%
_{i}\right) \vec{e}^{i}=\varkappa _{\alpha \beta }^{P}\left( \partial
_{q^{\delta }}\vec{w}_{P}\cdot \vec{e}_{i}\right) \vec{e}^{i},
\end{equation*}
are also vectors of the tangent vector space $\vec{g}_{\delta }$ of \textit{%
Riemannian} sub-space $q^{\alpha }$ of class $\hat{C}$%
\begin{equation*}
-\left( \partial _{q^{\delta }}\vec{N}_{\alpha \beta }\cdot \vec{e}%
_{i}\right) \vec{e}^{i}=\hat{\kappa}^{2}g_{\alpha \beta }e_{ki}\partial
_{q^{\delta }}x^{k}\vec{e}^{i}=\hat{\kappa}^{2}g_{\alpha \beta }\left( \vec{g%
}_{\delta }\cdot \vec{e}_{i}\right) \vec{e}^{i}=\hat{\kappa}^{2}g_{\alpha
\beta }\vec{g}_{\delta },
\end{equation*}
more exactly, 
\begin{equation*}
\varkappa _{\alpha \beta }^{P}t_{ij,P}\partial _{q^{\gamma }}x^{i}\partial
_{q^{\delta }}x^{j}=\varkappa _{\alpha \beta }^{P}\varkappa _{\gamma \delta
,P}=\hat{\kappa}^{2}g_{\alpha \beta }g_{\gamma \delta },
\end{equation*}
then, the normal curvature $\kappa _{n}$ of sub-space of curvature\textit{:} 
$q^{\alpha }$, of class $\hat{C}$, of the ambient space $x^{i}$ of class $C$%
, is determined by functional form 
\begin{equation*}
\left( \kappa _{n}\right) ^{2}=\frac{\left| \varkappa _{\alpha \beta
}^{P}\varkappa _{\delta ,P}^{\alpha }\right| }{\left| g_{\alpha \beta
}\right| }=\hat{\kappa}^{2M},
\end{equation*}
more exactly, 
\begin{equation*}
^{M}\sqrt{\left( \kappa _{n}\right) ^{2}}=\hat{\kappa}^{2}=\frac{\left(
t_{ij}^{P}t_{kl,P}-t_{il}^{P}t_{kj,P}\right) \partial _{q^{\alpha
}}x^{i}\partial _{q^{\beta }}x^{j}\partial _{q^{\gamma }}x^{l}\partial
_{q^{\delta }}x^{k}}{\left( e_{ij}e_{kl}-e_{il}e_{kj}\right) \partial
_{q^{\alpha }}x^{i}\partial _{q^{\beta }}x^{j}\partial _{q^{\gamma
}}x^{l}\partial _{q^{\delta }}x^{k}};
\end{equation*}
\begin{equation*}
^{M}\sqrt{\left( \kappa _{n}\right) ^{2}}=\hat{\kappa}^{2}=\frac{%
R_{ikjl}\partial _{q^{\alpha }}x^{i}\partial _{q^{\beta }}x^{j}\partial
_{q^{\gamma }}x^{l}\partial _{q^{\delta }}x^{k}}{\left(
e_{ij}e_{kl}-e_{il}e_{kj}\right) \partial _{q^{\alpha }}x^{i}\partial
_{q^{\beta }}x^{j}\partial _{q^{\gamma }}x^{l}\partial _{q^{\delta }}x^{k}}=
\end{equation*}
\begin{equation*}
=\frac{R_{\alpha \delta \beta \gamma }}{g_{\alpha \beta }g_{\gamma \delta
}-g_{\alpha \gamma }g_{\beta \delta }}=-\frac{R}{M\left( M-1\right) },
\end{equation*}
where the scalar invariant $R$\textit{:} $R=g^{\alpha \gamma }g^{\beta
\delta }R_{\alpha \delta \beta \gamma }$, is so-called invariant of
curvature (scalar curvature), \cite{An}.

It can be proved by an application of the result of \textit{Schur Theorem}, 
\cite{An}, that \textit{Riemannian} sub-spaces of curvature\textit{:} $%
q^{\alpha }$, of class $\hat{C}$, of the ambient space $x^{i}$ of class $C$,
whether they are principals or not, in this emphasized case are isotropic
spaces of the constant normal curvature $\kappa _{n}$.$\blacktriangledown $
\end{description}

\subsubsection{Example}

The curvature of \textit{Riemannian} spaces\textit{:} 
\begin{equation*}
ds^{2}=e^{\mu \left( \rho \right) }d\rho ^{2}+\rho ^{2}d\theta ^{2}+\sin
^{2}\theta d\varphi ^{2}+e^{\nu \left( \rho \right) }d\tau ^{2};
\end{equation*}
\begin{equation*}
\mu \left( \rho \right) =-\nu \left( \rho \right) =-\ln \left( 1-\frac{2m}{%
\rho }\right) ;\,m=const.
\end{equation*}
and 
\begin{equation*}
d\bar{s}^{2}=e^{\mu \left( \rho \right) }\left( d\rho ^{2}+\rho ^{2}d\theta
^{2}+\sin ^{2}\theta d\varphi ^{2}+d\tau ^{2}\right) ;
\end{equation*}
\begin{equation*}
\mu \left( \rho \right) =-\frac{2m}{\rho };\,m=const.,
\end{equation*}
with spherical symmetry.

Components of \textit{Christoffel's} symbols\textit{:} 
\begin{equation*}
\Gamma _{\alpha \beta ,\gamma }=\frac{1}{2}\left( \partial _{q^{\alpha
}}g_{\beta \gamma }+\partial _{q^{\beta }}g_{\alpha \gamma }-\partial
_{q^{\gamma }}g_{\alpha \beta }\right) ,
\end{equation*}
as well as of \textit{Reimann-Christoffel's} tensor of curvature: 
\begin{equation*}
R_{\alpha \beta \delta \gamma }=\partial _{q^{\delta }}\Gamma _{\beta \gamma
,\alpha }-\partial _{q^{\gamma }}\Gamma _{\beta \delta ,\alpha }+g^{\lambda
\sigma }\left( \Gamma _{\beta \delta ,\lambda }\Gamma _{\alpha \gamma
,\sigma }-\Gamma _{\beta \gamma ,\lambda }\Gamma _{\alpha \delta ,\sigma
}\right) ,
\end{equation*}
which are not identically zeros, for these sub-classes of general class of 
\textit{Riemannian} spaces with spherical symmetry, are the following forms%
\textit{:} 
\begin{equation*}
\Gamma _{11,1}=\bar{\Gamma}_{11,1}=\frac{1}{2}e^{\mu }\partial _{\rho }\mu
;\Gamma _{22,1}=-\Gamma _{12,2}=-\rho ;\bar{\Gamma}_{22,1}=-\bar{\Gamma}%
_{12,2}=-\frac{\rho e^{\mu }}{2}\left( 2+\rho \partial _{\rho }\mu \right) ;
\end{equation*}
\begin{equation*}
\Gamma _{33,1}=-\Gamma _{13,3}=-\rho \sin ^{2}\theta ;\bar{\Gamma}_{33,1}=-%
\bar{\Gamma}_{13,3}=-\frac{\rho \sin ^{2}\theta e^{\mu }}{2}\left( 2+\rho
\partial _{\rho }\mu \right) ;
\end{equation*}
\begin{equation*}
\Gamma _{33,2}=-\Gamma _{23,3}=-\rho \sin \theta \cos \theta ;\bar{\Gamma}%
_{33,2}=-\bar{\Gamma}_{23,3}=-\rho ^{2}e^{\mu }\sin \theta \cos \theta ;
\end{equation*}
\begin{equation*}
\Gamma _{44,1}=-\Gamma _{14,4}=-\frac{1}{2}e^{\nu }\partial _{\rho }\nu ;%
\bar{\Gamma}_{44,1}=-\bar{\Gamma}_{14,4}=-\frac{1}{2}e^{\mu }\partial _{\rho
}\mu ,
\end{equation*}
as well as 
\begin{equation*}
R_{1212}=\frac{1}{2}\rho \partial _{\rho }\mu ;\bar{R}_{1212}=\frac{1}{2}%
\rho e^{\mu }\left( \partial _{\rho }\mu +\partial _{\rho \rho }^{2}\mu
\right) ;
\end{equation*}
\begin{equation*}
R_{1313}=\frac{1}{2}\rho \sin ^{2}\theta \partial _{\rho }\mu ;\bar{R}%
_{1313}=\frac{1}{2}\rho e^{\mu }\sin ^{2}\theta \left( \partial _{\rho }\mu
+\rho \partial _{\rho \rho }^{2}\mu \right) ;
\end{equation*}
\begin{equation*}
R_{2323}=\rho ^{2}\sin ^{2}\theta (1-e^{-\mu });\bar{R}_{2323}=-\rho
^{3}e^{\mu }\partial _{\rho }\mu \sin ^{2}\theta \left( 1+\frac{\rho }{4}%
\partial _{\rho }\mu \right) ;
\end{equation*}
\begin{equation*}
R_{1414}=\frac{1}{2}e^{\nu }\left\{ \frac{1}{2}\left[ \partial _{\rho }\mu
\partial _{\rho }\nu -\left( \partial _{\rho }\nu \right) ^{2}\right]
-\partial _{\rho \rho }^{2}\nu \right\} ;\bar{R}_{1414}=-\frac{1}{2}e^{\mu
}\partial _{\rho \rho }^{2}\mu ;
\end{equation*}
\begin{equation*}
R_{2424}=-\frac{1}{2}\rho e^{\nu -\mu }\partial _{\rho }\nu ;\bar{R}_{2424}=-%
\frac{1}{2}\rho e^{\mu }\partial _{\rho }\mu \left( 1+\frac{\rho }{2}%
\partial _{\rho }\mu \right) .
\end{equation*}

By the functional relation (\ref{8}), it follows that\textit{:} 
\begin{equation*}
\kappa =\frac{1}{R_{1212}g}\sqrt{\left( R_{1212}\right)
^{2}R_{1313}R_{1414}R_{2323}R_{2424}}=\frac{1}{g}\sqrt{%
R_{1313}R_{1414}R_{2323}R_{2424}},
\end{equation*}
\begin{equation*}
\kappa =\frac{\left\{ \frac{1}{8}\rho ^{4}(1-e^{-\mu })e^{2\nu -\mu }\sin
^{4}\theta \partial _{\rho }\mu \partial _{\rho }\nu \left[ \partial _{\rho
\rho }^{2}\nu +\frac{1}{2}\left( \partial _{\rho }\nu \right) ^{2}-\frac{1}{2%
}\partial _{\rho }\mu \partial _{\rho }\nu \right] \right\} ^{\frac{1}{2}}}{%
\rho ^{4}e^{\mu +\nu }\sin ^{2}\theta }=2\frac{m^{2}}{\rho ^{6}}
\end{equation*}
and 
\begin{equation*}
\bar{\kappa}=\frac{1}{\bar{R}_{1212}\bar{g}}\sqrt{\left( \bar{R}%
_{1212}\right) ^{2}\bar{R}_{1313}\bar{R}_{1414}\bar{R}_{2323}\bar{R}_{2424}}=%
\frac{1}{\bar{g}}\sqrt{\bar{R}_{1313}\bar{R}_{1414}\bar{R}_{2323}\bar{R}%
_{2424}},
\end{equation*}
\begin{equation*}
\bar{\kappa}=\frac{\left[ \frac{1}{8}\rho ^{5}e^{4\mu }\sin ^{4}\theta
\left( \partial _{\rho }\mu \right) ^{2}\left( \partial _{\rho }\mu +\rho
\partial _{\rho \rho }^{2}\mu \right) \partial _{\rho \rho }^{2}\mu \left( 1+%
\frac{\rho }{4}\partial _{\rho }\mu \right) \left( 1+\frac{\rho }{2}\partial
_{\rho }\mu \right) \right] ^{\frac{1}{2}}}{\rho ^{4}e^{4\mu }\sin
^{2}\theta }=
\end{equation*}
\begin{equation*}
=2\frac{m^{2}}{\rho ^{6}}e^{\frac{4m}{\rho }}\sqrt{\left( 1+\frac{m}{2\rho }%
\right) \left( 1+\frac{m}{\rho }\right) }.\blacktriangledown
\end{equation*}

\section{CONCLUSION}

By functional form (\ref{8}), derived from general functional form (\ref{3})
generalizing the concept of the curvature of \textit{Riemannian} both one
and two dimensional spaces to the general concept of the curvature of 
\textit{Riemannian} vector spaces $q^{\alpha }$ of class $C$\textit{:} $%
C=N-M $, the concept itself of the curvature of \textit{Riemannian} vector
spaces of higher dimensions $\left( M\geq 2\right) $, directly related to
internal geometry of space, more exactly, to components of \textit{%
Reimann-Christoffel's} tensor of curvature $R_{\alpha \delta \beta \gamma }$.

The process of generalization of fundamental concepts of the differential
geometry, presented in this paper, gives the solid base to further
generalization other, whether they are fundamentals or not, concepts and
theorems of the differential geometry of a surface, and what may be the
subject of separated analysis.

The one of such concepts is that of \textit{Codazzi's} equations (formula)
of a surface, \cite{An} and \cite{Vi}, which can be, in the general case of 
\textit{Riemannian} space $q^{\alpha }$ of class $C$\textit{:} $C=N-M$,
obtained by projection of vector condition of integrability\textit{:} 
\begin{equation*}
\partial _{q^{\delta }q^{\beta }}^{2}\vec{g}_{\alpha }-\partial _{q^{\beta
}q^{\delta }}^{2}\vec{g}_{\alpha }=0,
\end{equation*}
of absolute differentials of the fundamental vectors $\vec{g}_{\alpha }$%
\textit{:}$\,d\vec{g}_{\alpha }=\partial _{q^{\beta }}\vec{g}_{\alpha
}dq^{\beta }$, onto the normal vector space of \textit{Riemannian} space $%
q^{\alpha }$ of class $C$.


\begin{thebibliography}{9}
\bibitem{An}  \-T. P. Andelic, \textit{Tensorial calculus (in Serbian), }%
Scientific book, Belgrade, 1980.

\bibitem{L-C}  T. Levi-Civita, \textit{The absolute differential calculus. }%
Blackie\&Son Ltd., Glasgow, 1954.

\bibitem{M-D}  D. S. Mitrinovic and D. Z. Dokovic, \textit{Polynomials and
matrices (in Serbian), }Civil book, Belgrade, 1986.

\bibitem{Pa}  V. Pauli, \textit{Theory of relativity (in Russian)}, Science,
Moscow, 1983.

\bibitem{Vi}  M. J. Vigodski, \textit{Differential geometry (in Russian)},
State department of publishing, Moscow, 1949.
\end{thebibliography}
\end{document}